\documentclass[a4paper,12pt]{amsart}
\usepackage{amsmath,amssymb,colordvi}
\usepackage{verbatim}
\newtheorem{theorem}{Theorem}

\newtheorem{lemma}[theorem]{Lemma}
\newtheorem{example}[theorem]{\it Example}
\newtheorem{proposition}[theorem]{Proposition}
\newtheorem{definition}[theorem]{Definition}

\newtheorem{remark}[theorem]{\it Remark}

\newcommand{\CaixaPreta}{\vrule Depth0pt height5pt width5pt}
\newcommand{\bgproof}{\noindent {\bf Proof.} \hspace{2mm}}
\newcommand{\edproof}{\hfill \CaixaPreta \vspace{3mm}}

\def\NN{\mathbb N}
\def\CC{\mathbb C}

\def\R{\mathbb R}
\def\RR{\mathbb R}

\def\S{\mathcal S}

\begin{document}

\title[Generalized Ces\`{a}ro operators on Sobolev spaces]
{ On the boundedness of generalized Ces\`{a}ro operators on
Sobolev spaces}

\author{Carlos Lizama}
\address{Departamento de Matem\'atica, Facultad de Ciencias, Universidad de Santiago de Chile, Casilla 307-Correo 2, Santiago-Chile.}

\email{carlos.lizama@usach.cl}

\thanks{C.Lizama and R. Ponce are partially supported by Project CONICYT-PIA ACT1112 Stochastic Analysis Research Network}

\author{Pedro J. Miana}
\address{Departamento de Matem\'aticas, Instituto Universitario de Matem\'aticas y Aplicaciones, Universidad de Zaragoza, 50009 Zaragoza, Spain.}
\email{pjmiana@unizar.es}

\author{Rodrigo Ponce}
\address{Instituto de Matem\'atica y F\'isica, Universidad de Talca, Casilla 747,
Talca, Chile.}
\email{rponce@inst-mat.utalca.cl}

\author{Luis S\'{a}nchez-Lajusticia}
\address{Departamento de Matem\'aticas, Instituto Universitario de Matem\'aticas y Aplicaciones, Universidad de Zaragoza,  50009 Zaragoza, Spain.}
\email{luiss@unizar.es}

\thanks{P.J. Miana and L. S\'{a}nchez-Lajusticia have been partially supported by Project MTM2010-16679, DGI-FEDER, of the MCYTS; Project E-64, D.G. Arag\'on, and JIUZ-2012-CIE-12, Universidad de Zaragoza, Spain.}




\begin{abstract}
For $\beta>0$ and $p\ge 1$, the generalized Ces\`{a}ro operator
$$
\mathcal{C}_\beta f(t):=\frac{\beta}{t^\beta}\int_0^t (t-s)^{\beta-1}f(s)ds
$$
and its companion operator $\mathcal{C}_\beta^*$ defined on Sobolev spaces $\mathcal{T}_p^{(\alpha)}(t^\alpha)$ and
 $\mathcal{T}_p^{(\alpha)}(\vert t\vert^\alpha)$ (where $\alpha\ge 0$ is the fractional order of derivation and are
 embedded in $L^p(\RR^+)$ and $L^p(\RR)$ respectively) are studied. We prove that if $p>1$, then $\mathcal{C}_\beta$ and $\mathcal{C}_\beta^*$ are
  bounded operators and commute on $\mathcal{T}_p^{(\alpha)}(t^\alpha)$ and $\mathcal{T}_p^{(\alpha)}(\vert t\vert^\alpha)$. We show explicitly the
   spectra $\sigma (\mathcal{C}_\beta)$ and  $\sigma (\mathcal{C}_\beta^*)$ and its operator norms (which depend on $p$). For $1< p\le 2$, we prove that $
\widehat{{\mathcal C}_\beta(f)}={\mathcal C}_\beta^*(\widehat{f})$ and $\widehat{{\mathcal C}_\beta^*(f)}={\mathcal C}_\beta(\widehat{f})$
where $\widehat{f}$ is the Fourier transform of a function $f\in L^p(\RR)$.
\end{abstract}

\date{}

\maketitle

\section{Introduction}

\setcounter{theorem}{0}
\setcounter{equation}{0}

Given $1\leq p<\infty,$ let $L^p(\mathbb{R}^+)$ be the set of Lebesgue $p$-integrable functions, that is, $f$ is a measurable function and
$$
||f||_p:=\left(\int_0^\infty |f(t)|^pdt\right)^{1/p}<\infty.
$$

The classical Hardy inequality (see \cite[p. 245]{Ha-Li-Po-64}) establishes that
$$
\left(\int_0^\infty \left|\frac{1}{t}\int_0^t f(s)ds\right|^pdt\right)^{1/p}\leq \frac{p}{p-1}||f||_p, \quad f\in L^p(\mathbb{R}^+),
$$
for $1<p<\infty$ and therefore the so-called Ces\`{a}ro transformation ${\mathcal C},$ defined by
\begin{equation}\label{eq1.1}
{\mathcal C}(f)(t)=\frac{1}{t}\int_0^t f(s)ds, \quad t> 0,
\end{equation}
is a bounded operator on $L^p(\mathbb{R}^+)$ with $||{\mathcal C}||\le \frac{p}{p-1}$ for $1<p<\infty$. In fact, it is also known that if  $\beta>0$
\begin{equation}\label{inicial}
\left(\int_0^\infty  \left|{\beta \over t^\beta} \int_0^t (t-s)^{\beta-1} f(s)ds\right\vert^pdt\right)^{1/p}\le {\Gamma(\beta +1)\Gamma(1-{1\over p})\over \Gamma(\beta +1-{1\over p})}\|f\|_p, \qquad f\in L^p(\RR^+),
\end{equation}
for $1<p<\infty$ and the constant ${\Gamma(\beta +1)\Gamma(1-{1\over p})\over \Gamma(\beta +1-{1\over p})}$ is optimal in this inequality, see \cite[Theorem 329]{Ha-Li-Po-64}. A closer (and dual) inequality is the following
\begin{equation}\label{dual}
\left(\int_{0}^{\infty} \!\! \,\left\vert\beta\int_{x}^{\infty} \!\! \frac{(t-x)^{\beta-1}}{t^\alpha} f(t)dt\right\vert^{p}dx\right)^{1\over p} \leq \frac{\Gamma(\alpha+1)\Gamma\left(\frac{1}{p}\right)}{\Gamma\left(\alpha+\frac{1}{p}\right)}\Vert f\Vert_p.
\end{equation}
Also the constant $\frac{\Gamma(\alpha+1)\Gamma\left(\frac{1}{p}\right)}{\Gamma\left(\alpha+\frac{1}{p}\right)} $ is optimal in the  above inequality (\cite[Theorem 329, p.245]{Ha-Li-Po-64}).

Note that inequalities (\ref{inicial}) and (\ref{dual}) show that operators $\mathcal{C}_\beta$, $\mathcal{C}_\beta^*$ where
$$
\mathcal{C}_\beta f(t):=\frac{\beta}{t^\beta}\int_0^t (t-s)^{\beta-1}f(s)ds, \qquad \mathcal{C}_\beta^* f(s):=\beta\int_s^\infty \frac{(t-s)^{\beta-1}}{t^\beta}f(t)dt,
$$
define bounded operators on $L^p(\RR^+)$, $\mathcal{C}_1=\mathcal{C}$ and $\mathcal{C}_1^*=\mathcal{C}^*$. By Fubini theorem, the dual operator of $\mathcal{C}_\beta$ on $L^p(\RR^+)$ is $\mathcal{C}_\beta^*$ on $L^{p'}(\RR^+)$, i.e,
$$
\int_0^\infty \mathcal{C}_\beta f(t)g(t)dt=\int_0^\infty f(s)\mathcal{C}_\beta^{*} g(s)ds, \qquad f \in L^p(\RR^+), \quad g\in L^{p'}(\RR^+),
$$
where $1<p, p'<\infty$ and ${1\over p}+{1\over p'}=1 $. See other properties about some of these operators in \cite{Br-Ha-Sh-65, [Bo], Mo99}.

Recently, A. Arvanitidis and A. Siskakis (\cite{Ar-Si-11}) showed that the half-plane versions of Ces\`{a}ro operators on the Hardy space ${\mathcal H}_p(\mathbb{U}),$ defined on $\mathbb{U}:=\{z\in \mathbb{C}:{\rm Im}(z)>0\}$ by
\begin{equation}\label{eq1.2}
{C}(F)(z):=\frac{1}{z}\int_0^z F(s)ds, \quad {C}^*(F)(z):=\int_z^\infty {F(s)\over s} ds, \quad F\in H^p(\mathbb{U}),
\end{equation}
define bounded operators on ${\mathcal H}_p(\mathbb{U})$ when $p>1$. Both operators ${C}$ and ${C}^*$ can be obtained as resolvent operators of generators of some appropriate strongly continuous $C_0$-semigroups on ${\mathcal H}_p(\mathbb{U}).$

Similarly, W. Arendt and B. de Pagter (\cite{Ar-Pa-02}) studied the Ces\`{a}ro operator (\ref{eq1.1}) defined in an interpolation space $E$ of $(L^1,L^\infty)$ on $\mathbb{R}^+.$ When $E=L^p(\mathbb{R}^+),$ the authors obtain a representation of $\mathcal{C}$ in terms of an appropriate resolvent operator, see \cite[Corollaries 2.2, 4.3]{Ar-Pa-02}.

In \cite{Ga-Mi-06},  Sobolev subspaces
$\mathcal{T}^{(\alpha)}_1(t^\alpha)$ and
$\mathcal{T}^{(\alpha)}_1(\vert t\vert^\alpha)$ (contained in
$L^1(\mathbb{R}^+)$ and $L^1(\mathbb{R})$ respectively and where
$\alpha\ge 0$ is the fractional order of derivation) were
introduced. In fact, these subspaces are sub-algebras for the
convolution products given by
\begin{equation}\label{convos-intro}
f\ast g(t)=\int_0^tf(t-s)g(s)ds, \quad t\geq0,
\end{equation}
and
\begin{equation}\label{convos2-intro}
f\ast g(t)=\int_{-\infty}^{\infty} f(t-s)g(s)ds, \qquad t\in \R,
\end{equation}
respectively. These algebras  are canonical to define some algebra homomorphisms (defined by integral representations) into $\mathcal{B}(X),$ the set of all linear and bounded operators on a Banach space $X$. See further details in \cite{Ga-Mi-06}.

Further, in \cite{Ro-08}   Sobolev subspaces
$\mathcal{T}^{(\alpha)}_p(t^\alpha)$ contained in Lebesgue spaces
$L^p(\RR^+)$ ($p\ge 1$) were introduced and studied in detail.
Some remarkable results were proved (see Proposition \ref{juanjo}
below). In particular, the subspace
$\mathcal{T}^{(\alpha)}_p(t^\alpha)$ is a module for the algebra
$\mathcal{T}^{(\alpha)}_1(t^\alpha)$ for the convolution product
$\ast$ given by (\ref{convos-intro}).

Hence, it is natural to ask in what extent the boundedness
property of the operators $\mathcal{C}_\beta$ and
$\mathcal{C}_\beta^*$  remain valid in the above described Sobolev
spaces.

The main aim of this paper is to study boundedness, representation
and spectral properties for the generalized Ces\`{a}ro operators
$\mathcal{C}_\beta$ and $\mathcal{C}_\beta^*$ on Sobolev subspaces
of fractional order $\alpha \ge 0$ embedded in $L^p(\RR^+)$ and
$L^p(\RR)$ (which are denoted by
$\mathcal{T}^{(\alpha)}_p(t^\alpha)$ and
$\mathcal{T}^{(\alpha)}_p(\vert t\vert^\alpha)$ respectively).

The outline of the paper is as follows: In the second section we remind some basic properties of the Sobolev spaces $\mathcal{T}_p^{(\alpha)}(t^\alpha)$ (where $\mathcal{T}_p^{(\alpha)}(t^\alpha)\hookrightarrow L^p(\mathbb{R}^+)) $. We also prove new results, see for example
Proposition \ref{continu}. The main tool of this section (and in the rest of the paper) are the $C_0$-group of isometries on
$\mathcal{T}_p^{(\alpha)}(t^\alpha)$, $(T_{t, p})_{t\in \RR}$ given by
$$
T_{t,p} f(s):=e^{-\frac{t}{p}}f(e^{-t}s), \qquad f \in \mathcal{T}_p^{(\alpha)}(t^\alpha).
$$
In the Theorem \ref{semigroup} it is identified its infinitesimal
generator and, its spectrum, in Proposition \ref{spectrum}. We
note that this strategy has been pursued by other authors. We
mention here \cite{Ar-Pa-02, Ar-Si-11, CC, [Si]}.

In the third section, we study the generalized Ces\`{a}ro
operators $\mathcal{C}_\beta$ and $\mathcal{C}_\beta^*$ defined on
Sobolev spaces $\mathcal{T}_p^{(\alpha)}(t^\alpha)$. We first show
that both operators are bounded operators and commute for $p>1$.
In fact, we have
$$
||\mathcal{C}_{\beta}|| = \frac{\Gamma(\beta +1)\Gamma(1/p')}{\Gamma(\beta +
1/p')}; \quad ||\mathcal{C}_{\beta}^* || = \frac{\Gamma(\beta
+1)\Gamma(1/p)}{\Gamma(\beta + 1/p)},
$$
for $\alpha \geq 0, p>1, \beta >0, \,\, 1/p +1/p' =1.$ It is
remarkable that the composition
$\mathcal{C}_\alpha\mathcal{C}_\beta^\ast$ may be described
explicitly involving the Gaussian hypergeometric function $_2F_1$
(see Theorem \ref{conmutacion}) as follows:
\begin{eqnarray*}
(\mathcal{C}_\alpha\mathcal{C}_\beta^\ast)f(t)&=&\alpha\int^t_0f(r){1\over
t-r}\left({t-r\over t}\right)^{\alpha+\beta} \,_2F_1(\alpha+\beta,
\beta; \beta+1; {r\over t})dr\cr\qquad &\quad&\qquad + \,\,
\beta\int_t^\infty f(r){1\over r-t}\left({r-t\over
t}\right)^{\alpha+\beta} \,_2F_1(\alpha+\beta,  \alpha; \alpha+1;
{t\over r})dr,
\end{eqnarray*}
for $\alpha, \beta >0.$

 Using the description of $\mathcal{C}_\beta$ and
$\mathcal{C}_\beta^*$ in terms of the $C_0$-semigroups (Theorem
\ref{lemma1.1} and Theorem \ref{theorem3.6}), we are able to
determine the spectra, $\sigma (\mathcal{C}_\beta)$ and $\sigma
(\mathcal{C}_\beta^*)$ (Theorem \ref{spec} and \ref{spectrr}) as:
$$
\sigma(\mathcal{C}_\beta)=
\Gamma(\beta+1)\overline{\left\{{\Gamma({1\over p'}+it)\over
\Gamma(\beta+{1\over p'}+it)} \ : \ t\in \mathbb{R}\right\}};
$$
and
$$
\sigma(\mathcal{C}_\beta^*)=
\Gamma(\beta+1)\overline{\left\{{\Gamma({1\over p}+it)\over
\Gamma(\beta+{1\over p}+it)} \ : \ t\in \mathbb{R}\right\}},
$$
where $1/p +1/p'=1.$ In particular, the operators $\mathcal{C}_1 $
and $\mathcal{C}_1^* $ can be obtained as the resolvent operator
of appropriate $C_0$-semigroups, namely $(T_{t, p})_{t\ge 0}$ and
$(T_{-t, p})_{t\ge 0},$ respectively.

We remark that in case $\beta=1$ we obtain:
\begin{equation*}
\sigma(\mathcal{C}_1^*)=\left\{w\in\mathbb{C} \ : \
\left|w-\frac{p}{2}\right|=\frac{p}{2}\right\}.
\end{equation*}
This  gives a proof of a conjecture posed by F. M\'{o}ricz on $L^p(\RR^+)$
\cite[Section 2]{Mo99} and  new proofs of some results given in \cite{Br-Ha-Sh-65,[Bo]}.

In Section 4, we introduce and give some basic properties of the
Sobolev spaces $\mathcal{T}_p^{(\alpha)}(\vert t\vert^\alpha)$
(here $\mathcal{T}_p^{(\alpha)}(\vert
t\vert^\alpha)\hookrightarrow L^p(\mathbb{R})) $. We also prove
that the space $\mathcal{T}_p^{(\alpha)}(\vert t\vert^\alpha)$ is
a module for the algebra $\mathcal{T}_1^{(\alpha)}(\vert
t\vert^\alpha)$ and the $\ast$-convolution product given by
(\ref{convos2-intro}). Moreover, the following interesting
inequality holds:
$$
|||f\ast g|||_{\alpha, p}\le C_{\alpha, p} |||f|||_{\alpha, p}|||
g|||_{\alpha, 1}, \qquad f\in
\mathcal{T}_p^{(\alpha)}(|t|^\alpha), \quad g \in
\mathcal{T}_1^{(\alpha)}(|t|^\alpha).
$$
In section 5, we study boundedness, representation and spectral
properties of generalized C\`esaro operators on $\mathbb{R}.$
Again, it is relevant to mention that the $C_0$-group of
isometries on $\mathcal{T}_p^{(\alpha)}(\vert t\vert^\alpha)$,
$(T_{t, p})_{t\in \RR}$ given by
$$
T_{t,p} f(s):=e^{-\frac{t}{p}}f(e^{-t}s), \qquad f \in \mathcal{T}_p^{(\alpha)}(\vert t\vert^\alpha),
$$
(Theorem \ref{groups}) is the main tool to prove the main results
in this section. The generalized Ces\`{a}ro operators
$\mathcal{C}_\beta$ and $\mathcal{C}_\beta^*$ defined on Sobolev
spaces $\mathcal{T}_p^{(\alpha)}(\vert t\vert^\alpha)$ are
described in terms of the $C_0$-group of isometries $(T_{t,
p})_{t\in \RR}$. Similar results shown in the case
$\mathcal{T}_p^{(\alpha)}( t^\alpha)$ hold in this case, see
Theorem \ref{lemma2.1} and \ref{lemma2.2} below.

In the last section we show
 that $\widehat{{\mathcal C}_\beta(f)}={\mathcal C}_\beta^*(\widehat{f})$ and $\widehat{{\mathcal C}_\beta^*(f)}={\mathcal C}_\beta(\widehat{f})$
  where $\widehat{f}$ is the Fourier transform of a function $f\in L^p(\RR)$ and $1< p\le 2$, see Theorem \ref{conmutan}. We notice that
   our studies in this section extends and complement the main result in \cite{Mo02}.

\section{Composition groups on Sobolev spaces defined on $\mathbb{R}^+.$}

\setcounter{theorem}{0}
\setcounter{equation}{0}

Let $\mathcal{D}_+$ be the class of $C^\infty$-functions with compact support on $[0,\infty)$ and $\mathcal{S}_+$ the Schwartz class on $[0,\infty)$. For a function $f\in \mathcal{S}_+$ and $\alpha>0$, the \textit{Weyl fractional integral} of order $\alpha$, $W^{-\alpha}_+f$, is defined by
$$
W^{-\alpha}_+f(t):=\frac{1}{\Gamma(\alpha)}\int_t^\infty (s-t)^{\alpha-1}f(s)ds, \qquad t\in\mathbb{R}^+.
$$
The \textit{Weyl fractional derivative} $W^\alpha_+f$ of order $\alpha$ is defined by
$$
W^\alpha_+ f(t):=(-1)^n\frac{d^n}{dt^n}W_+^{-(n-\alpha)}f(t),\quad t\in \mathbb{R}^+
$$
where $n=[\alpha]+1,$ and $[\alpha]$ denotes the integer part of $\alpha.$ It is proved that $W_+^{\alpha+\beta}=W_+^{\alpha}(W_+^{\beta})$ for any $\alpha, \beta\in\mathbb{R},$ where $W_+^0=Id$ is the identity operator and $(-1)^nW_+^n=\frac{d^n}{dt^n}$ holds with $n\in\mathbb{N}$, see more details in \cite{Mi-Ro} and \cite{Sa-Ki-Ma-93}.


Take $ \lambda>0$ and $f_\lambda$ defined by $f_\lambda(r):=f(\lambda r)$ for $r>0$ and $f \in \mathcal{S}_+$. It is direct to check that
\begin{equation}\label{escala}
W^\alpha_+f_\lambda = \lambda^\alpha (W^\alpha_+f)_\lambda, \qquad f \in \mathcal{S}_+,
\end{equation}
for $\alpha \in \RR$.





Now we introduce a family of subspaces $\mathcal{T}_p^{(\alpha)}(t^\alpha)$ which are contained in $L^p(\RR^+)$.

\begin{definition}
For $\alpha>0$ let be the Banach space $\mathcal{T}_p^{(\alpha)}(t^\alpha)$ defined as the completion of the Schwartz class  $\mathcal{S}_+$ in the norm
\begin{equation*}
||f||_{\alpha, p}:={1\over \Gamma(\alpha+1)}\left(\int_0^\infty |W^\alpha_+ f(t)|^pt^{\alpha p}dt\right)^{\frac{1}{p}}.
\end{equation*}
\end{definition}

We understand that  $\mathcal{T}_p^{(0)}(t^0)= L^p(\mathbb{R}^+)$ and $||\quad||_{0, p}=||\quad||_{p}$.  The case $p=1$ and $\alpha \in \NN$ where introduced in \cite{AK} and for $\alpha >0$ in \cite{Ga-Mi-06}.

In the next proposition we collect some  results about these family of spaces $\mathcal{T}_p^{(\alpha)}(t^\alpha)$ which we  may be found in \cite{Ro-08}.

\begin{proposition}\label{juanjo}
Take $p\geq 1$ and $\beta>\alpha>0$. Then
\begin{itemize}
\item[(i)] $\mathcal{T}_p^{(\beta)}(t^\beta)\hookrightarrow \mathcal{T}_p^{(\alpha)}(t^\alpha)\hookrightarrow  L^p(\mathbb{R}^+) $.

\item[(ii)] $\mathcal{T}_p^{(\alpha)}(t^\alpha)\ast \mathcal{T}_1^{(\alpha)}(t^\alpha) \hookrightarrow \mathcal{T}_p^{(\alpha)}(t^\alpha)$ for $1\le p <\infty$, where
\begin{equation}\label{convos}
f\ast g(t)=\int_0^tf(t-s)g(s)ds, \quad t\ge 0, \qquad f\in \mathcal{T}_p^{(\alpha)}(t^\alpha), \quad g\in\mathcal{T}_1^{(\alpha)}(t^\alpha).
\end{equation}

\item[(iii)] The operator $D^\alpha_+: \mathcal{T}_p^{(\alpha)}(t^\alpha)\to L^p(\RR^+)$  defined by
$$
f\mapsto D^\alpha_+ f(t)={1\over \Gamma(\alpha+1)}t^\alpha W^\alpha_+f(t), \qquad t\ge 0, \quad f \in \mathcal{T}_p^{(\alpha)}(t^\alpha).
$$
is an isometry.

\item[(iv)] If $p>1$ and $p'$ satisfies $\frac{1}{p}+\frac{1}{p'}=1$, then the dual of $\mathcal{T}_p^{(\alpha)}(t^\alpha)$ is $\mathcal{T}_{p'}^{(\alpha)}(t^\alpha)$, where the duality is given by $$\langle f,g \rangle_{\alpha} ={1\over \Gamma(\alpha+1)^2}\int_0^\infty W^\alpha_+ f(t)W^\alpha_+ g(t)t^{2\alpha}dt,$$ for $f\in \mathcal{T}_p^{(\alpha)}(t^\alpha)$, $g\in \mathcal{T}_{p'}^{(\alpha)}(t^\alpha)$.
\end{itemize}
\end{proposition}

Note that, in fact,
\begin{equation}
\Vert f\Vert_{\alpha,p}=\Vert D^\alpha_+ f\Vert_p, \qquad  \langle f,g \rangle_{\alpha}= \langle D^\alpha_+ f,D^\alpha_+ g \rangle_{0},
\end{equation}
for $f\in \mathcal{T}_p^{(\alpha)}(t^\alpha)$ and  $g\in \mathcal{T}_{p'}^{(\alpha)}(t^\alpha)$ with $\frac{1}{p}+\frac{1}{p'}=1$.


In the next lemma, we consider some  functions which belong (or not) to $\mathcal{T}_p^{(\alpha)}(t^\alpha)$ for $p\ge 1$.
\begin{lemma}\label{lemma2.3}
If $\alpha,a>0$ and $p\geq 1,$ then

(i) $t^{\beta}\not \in \mathcal{T}_p^{(\alpha)}(t^\alpha)$ for $\beta\in\mathbb{C}.$

(ii) $(a+t)^{-\beta}\in \mathcal{T}_p^{(\alpha)}(t^\alpha)$ for ${\Re}\beta>1/p.$
\end{lemma}

\bgproof
(i) It suffices to note that $t^{\beta}$ does not belong to $L^p(\mathbb{R}^+)$.

(ii) For $0<{\Re}\gamma<{\Re}\delta$ and $a>0$ it is well know that $W^{-\gamma}_+(a+t)^{-\delta}=\frac{\Gamma(\delta-\gamma)}{\Gamma(\delta)}(t+a)^{\gamma-\delta}$, see for example \cite[p. 201]{Er-Ma-Ob-Tr-54}. With this formula, it is easy to check that
$$
W^\alpha_+(a+t)^{-\beta}=\frac{\Gamma(\alpha+\beta)}{\Gamma(\beta)} (t+a)^{-(\alpha+\beta)}.
$$
Thus for $f(t):=(a+t)^{-\beta}$ we obtain
\begin{eqnarray*}
||f||_{\alpha,p}^p&=&{1\over \Gamma(\alpha+1)^p}\int_0^\infty |W^\alpha_+ f(t)|^p t^{\alpha p} dt=\left(\frac{\Gamma(\alpha+\beta)}{\Gamma(\alpha)\Gamma(\beta)} \right)^p\int_0^\infty \frac{t^{\alpha p}}{|(t+a)^{(\alpha+\beta)p}|} dt\\
&\leq&\left(\frac{\Gamma(\alpha+\beta)}{\Gamma(\alpha)\Gamma(\beta)} \right)^p\int_0^\infty \frac{1}{(t+a)^{p\Re\beta }} dt<\infty,
\end{eqnarray*}
and we conclude the proof.
\edproof

Given $f\in \mathcal{T}_p^{(\alpha)}(t^\alpha)$, as next result shows, we obtain that the function $f\in C(\R^+)$ for $p, \alpha\ge 1$.

\begin{proposition} \label{continu} Take $p, \alpha\ge 1$ and $f\in \mathcal{T}_p^{(\alpha)}(t^\alpha)$. Then $f\in C(\R^+)$, $\lim_{t\to \infty}f(t)=0$ and
$$
\sup_{t>0} t^{p} \vert f(t)\vert \le C_{\alpha,p}\Vert f\Vert_{\alpha, p},  \qquad f\in \mathcal{T}_p^{(\alpha)}(t^\alpha),
$$
where $C_{\alpha,p}$ is independent of $f$.
\end{proposition}

\bgproof
By Proposition \ref{juanjo} (i), it is enough to check for $\alpha=1$. Take $t > s > 0,$ and we get that
$$
\vert f(t)-f(s)\vert\le \int^t_s\vert f'(u)\vert du\le {1\over s}\int^t_s\vert f'(u)\vert u du.
$$
For $p=1$, it is clear that $f$ is continuous and for $p>1$, we apply the H\"{o}lder inequality to obtain
$$
\vert f(t)-f(s)\vert \le \Vert f\Vert_{1, p}\left(t-s\right)^{1\over p'}, \qquad {1\over p}+{1\over p'}=1.
$$
Then $f$ is continuous in $\R^+$.  Take $f\in \mathcal{T}_1^{(\alpha)}(t^\alpha)$, and
$$
\vert f(t)\vert \le \int_t^\infty\vert f'(u)\vert du\le {1\over t}\int_t^\infty u\vert f'(u)\vert du, \qquad t>0,
$$
and we conclude that  $\lim_{t\to \infty}f(t)=0$. Similarly take  $f\in \mathcal{T}_p^{(\alpha)}(t^\alpha)$ with $1<p<\infty$. Then we have that
$$
\vert f(t)\vert\le \int_t^\infty \vert f'(u)\vert du\le\left(\int_t^\infty u^p\vert f'(u)\vert^p du\right)^{1\over p}\left(\int_t^\infty {1\over u^{p'}}du\right)^{1\over p'}\le \left({1\over p't^{p'-1}}\right)^{1\over p'}\Vert f\Vert_{1, p}
$$
where we conclude that  $
\sup_{t>0} t^{p} \vert f(t)\vert \le \left({1\over p'}\right)^{1\over p'}\Vert f\Vert_{1, p}$ and the proof is finished.
\edproof

The following is the main result of this section. It will be the
key in the study of spectral properties of the  generalized
Ces\`{a}ro operators $\mathcal{C}_\beta$ and $\mathcal{C}_\beta^*$
defined on Sobolev spaces.

\begin{theorem} \label{semigroup} For $1\le p$ and $\alpha\ge 0$, the family of operators  $(T_{t,p})_{t\in \RR}$ defined by
$$
T_{t,p} f(s):=e^{-\frac{t}{p}}f(e^{-t}s), \qquad f \in \mathcal{T}_p^{(\alpha)}(t^\alpha),
$$
is a $C_0$-group of isometries on $\mathcal{T}_p^{(\alpha)}(t^\alpha)$ whose  infinitesimal generator $\Lambda$ is given by
$$
(\Lambda f)(s):=-sf'(s)-\frac{1}{p}f(s)
$$
with domain $D(\Lambda)= \mathcal{T}_p^{(\alpha+1)}(t^{\alpha+1}).$
\end{theorem}

\bgproof
We check that the operators $(T_{t,p})_{t\in \RR}$ are isometries:
\begin{eqnarray*}
||T_{t,p} f||_{\alpha,p}^p&=&{1\over \Gamma(\alpha+1)^p}\int_0^\infty |W^\alpha_+ T_{t,p} f(s)|^ps^{\alpha p}ds
={e^{-t}\over \Gamma(\alpha+1)^p}\int_0^\infty |W^\alpha_+ f(e^{-t}s)|^ps^{\alpha p}ds\\
&=&{e^{-t}\over \Gamma(\alpha+1)^p}\int_0^\infty e^t|e^{-\alpha t}(W^\alpha_+ f)(u)|^p(e^{\alpha t}u^\alpha)^pdu
=||f||_{\alpha,p}^p,
\end{eqnarray*}
where we have applied the equality (\ref{escala}).

Using some known properties for fractional derivative (\cite[p. 96]{Sa-Ki-Ma-93}) it can be shown that the family of operatos $(T_{t,p})_{t\in \RR}$ are strongly continuous, see similar ideas in \cite[Proposition 2.1]{Ar-Si-11} and \cite[Section 2]{Ar-Pa-02}. It is straightforward to check that the family $(T_{t,p})_{t\in \RR}$ is a group of operators.

On $\mathcal{T}_p^{(\alpha)}(t^\alpha)$ define $\{S_t\}_{t\geq 0}$ by $S_t(f)(s):=f(e^{-t}s).$ Then, an easy computation shows that the generator $A$ of $\{S_t\}_{t\geq 0}$ with domain $\{f\in \mathcal{T}_p^{(\alpha)}(t^\alpha) \ : \ tf'\in \mathcal{T}_p^{(\alpha)}(t^\alpha)\}$ is given by $Af(s)=-sf'(s).$ Therefore, the rescaled semigroup $(T_{t,p})_{t\geq 0}$ has domain $\{f\in \mathcal{T}_p^{(\alpha)}(t^\alpha) \ : \ tf'\in \mathcal{T}_p^{(\alpha)}(t^\alpha)\}$ and his generator is $(\Lambda f)(s)=-sf'(s)-\frac{1}{p}f(s).$ See \cite[p. 60]{En-Na-00} for more details.

Finally, we prove that $D(\Lambda)=\mathcal{T}_p^{(\alpha+1)}(t^{\alpha+1}).$ In fact, let $f\in \mathcal{T}_p^{(\alpha+1)}(t^{\alpha+1})$ be given. Since $\mathcal{T}_p^{(\alpha+1)}(t^{\alpha+1})\hookrightarrow \mathcal{T}_p^{(\alpha)}(t^\alpha)$, we have $f\in
\mathcal{T}_p^{(\alpha)}(t^\alpha)$. From \cite[p. 246]{Mi-Ro} it is easy to show that $W^\alpha_+ (tf'(t))=\alpha W^\alpha_+
f(t)+tW^{\alpha+1}_+ f(t).$ Thus, $tf'\in\mathcal{T}_p^{(\alpha)}(t^\alpha)$ and therefore $f\in D(\Lambda)$. Conversely, if $f\in D(\Lambda)$, then $f\in \mathcal{T}_p^{(\alpha)}(t^\alpha)$ and $tf'\in \mathcal{T}_p^{(\alpha)}(t^\alpha)$. The same above identity, implies that $t^{\alpha+1}W^{\alpha+1}_+f(t)=t^\alpha W^{\alpha}_+ (tf'(t))-\alpha t^\alpha W^\alpha_+ f(t)$, and therefore
$f\in\mathcal{T}_p^{(\alpha+1)}(t^{\alpha+1})$.
\edproof

The proof of the following result is inspired in \cite[Proposition 2.3]{Ar-Si-11}. We denote by $\sigma(\Lambda)$ the usual spectrum of the operator $\Lambda$ and by $\sigma_p(\Lambda)$ the point spectrum of the operator $\Lambda$.

\begin{proposition} \label{spectrum} For $1\leq p<\infty$ we have

$(i)$ $\sigma_p(\Lambda)=\emptyset;$

$(ii)$ $\sigma(\Lambda)=i\mathbb{R}.$
\end{proposition}

\bgproof
$(i)$ Let $\lambda\in\mathbb{C}$ and $f\in \mathcal{T}_p^{(\alpha)}(t^\alpha)$  such that $\Lambda(f)=\lambda f$. Then, $f$ is solution of the differential equation
$$
sf'(s)+(\lambda+\frac{1}{p})f(s)=0.
$$
The nonzero solutions to this equation have the form $f(t)=ct^{-(\lambda+1/p)}$ with $c\neq 0.$ But by Lemma \ref{lemma2.3}, these solutions are not in $\mathcal{T}_p^{(\alpha)}(t^\alpha).$ Therefore $\sigma_p(\Lambda)=\emptyset.$

$(ii)$ Since each $T_{t,p}$ is an invertible isometry its spectrum satisfies
$$
\sigma(T_{t,p})\subseteq \{z\in\mathbb{C}:|z|=1\}.
$$
By the spectral mapping theorem (see Theorem \cite[IV.3.6]{En-Na-00}), we have that
$$
e^{t\sigma(\Lambda)}\subseteq \sigma (T_{t,p}).
$$
Therefore, if $w\in \sigma(\Lambda),$ then $e^{tw}\in  \{z\in\mathbb{C}:|z|=1\}.$ Thus, we obtain that $\sigma(\Lambda)\subseteq i\mathbb{R}.$

Conversely, let $\mu\in i\mathbb{R}$ and assume that $\mu\in\rho(\Lambda).$ Let $\lambda=\mu+\frac{1}{p}.$ By Lemma \ref{lemma2.3} the function $f$ defined by $f(t):=(1+t)^{-\lambda-1}\in \mathcal{T}_p^{(\alpha)}(t^\alpha).$ Since $R(\mu,\Lambda)$ is a bounded operator, the function $g(t):=R(\mu,\Lambda)f(t)$ belongs to $\mathcal{T}_p^{(\alpha)}(t^\alpha).$ Therefore, $g$ is solution of equation
$$
\lambda g(t)+tg'(t)=f(t).
$$
An easy computation shows that the solution of this equation is $G(t):=ct^{-\lambda}+\lambda^{-1}(1+t)^{-\lambda},$  where $c$ is a constant. However, as in Lemma \ref{lemma2.3} one can check that $G\not\in \mathcal{T}_p^{(\alpha)}(t^\alpha).$ Therefore, $\mu\in \sigma(\Lambda).$
\edproof

Now, consider the negative part $\{T_{-t,p}, t\ge 0\}$ of the group $\{T_{t,p}\}_{t\in \mathbb{R}}$: that is, for $f\in \mathcal{T}_p^{(\alpha)}(t^\alpha)$,
$$
T_{-t,p}f(s)=e^{\frac{t}{p}}f(e^ts), \,\, t\geq 0.
$$
Obviously, $\{T_{-t,p}\}_{t\geq 0}$ is a $C_0$-semigroup on $\mathcal{T}_p^{(\alpha)}(t^\alpha)$  of isometries whose generator is $-\Lambda$.

We finish this section, establishing the relationship between the semigroups $\{T_{t,p}\}_{t\geq 0}$ and $\{T_{-t, p'}\}_{t\geq 0}$ with ${1\over p}+ {1\over p'}=1.$

\begin{proposition}\label{dualsemi}
The semigroups $\{T_{t,p}\}_{t\geq 0}$ and $\{T_{-t,p'}\}_{t\geq 0}$ are dual operators of each other acting on $\mathcal{T}_p^{(\alpha)}(t^\alpha)$ and $\mathcal{T}_{p'}^{(\alpha)}(t^\alpha)$ with  ${1\over p}+ {1\over p'}=1.$
\end{proposition}

\bgproof
By the Proposition \ref{juanjo} (iv) the  dual of $\mathcal{T}_p^{(\alpha)}(t^\alpha)$ is  $\mathcal{T}_{p'}^{(\alpha)}(t^\alpha)$ (with ${1\over p}+ {1\over p'}=1.)$ Take $f\in \mathcal{T}_p^{(\alpha)}(t^\alpha)$ and $g\in  \mathcal{T}_{p'}^{(\alpha)}(t^\alpha)$. Then
\begin{eqnarray*}
\langle T_{t,p} f,g\rangle_\alpha &=&{1\over \Gamma(\alpha +1)^2}\int_0^\infty (W^\alpha_+ T_{t,p} f)(s)W^\alpha_+ g(s)s^{2\alpha}{ds}\\
&=&{1\over \Gamma(\alpha +1)^2}\int_0^\infty e^{-\frac{t}{p}}e^{-t\alpha}(W^\alpha_+ f)(e^{-t}s)W^\alpha_+ g(s)s^{2\alpha}{ds}\\
&=&{1\over \Gamma(\alpha +1)^2}\int_0^\infty e^{-\frac{t}{p}}e^t(W^\alpha_+ f)(u)e^{\alpha t}(W^\alpha_+ g)(ue^t)u^{2\alpha}{du}\\
&=&{1\over \Gamma(\alpha +1)^2}\int_0^\infty (W^\alpha_+ f)(u)(W^\alpha_+ T_{-t,p'}g)(u)u^{2\alpha}{du}=\langle f, T_{-t,p'} g\rangle_\alpha,
\end{eqnarray*}
where we change the variable and we conclude the result.
\edproof

\section{Generalized Ces\`{a}ro operators on Sobolev spaces defined on $\mathbb{R}^+.$}

\setcounter{theorem}{0}
\setcounter{equation}{0}

For $\beta>0$ the generalized Ces\`{a}ro operator  on $\mathcal{T}_p^{(\alpha)}(t^\alpha)$ is
defined by
$$
\mathcal{C}_\beta f(t):=\frac{\beta}{t^\beta}\int_0^t (t-s)^{\beta-1}f(s)ds=\beta\int_0^1(1-r)^{\beta-1}f(tr)dr, \quad \,\,t> 0.
$$

Defining the function
$$
g_{\beta}(t)= \frac{t^{\beta-1}}{\Gamma(\beta)}, \quad t>0,
$$
we obtain the also  equivalent formulation of the generalized Ces\`{a}ro operator in terms of  finite convolution as follows:
$$
\mathcal{C}_\beta f(t):=\frac{1}{g_{\beta + 1}(t)}\int_0^t g_{\beta}(t-s)f(s)ds, \quad t>0.
$$
We remark that for certain classes of vector-valued functions $f$, the asymptotic behavior as $t \to \infty$ of $\mathcal{C}_\beta f(t)$ in the above representation has been studied in \cite{Li-Pr-03}.

Note that we may calculate ${\mathcal C}_\beta (f)$ for some particular functions:

\begin{example}
{\rm (i) Functions $g_\gamma$ are eigenfunctions of ${\mathcal C}_\beta$ with eigenvalue ${\Gamma(\beta+1) \Gamma(\gamma)\over \Gamma(\beta+\gamma)}$:
$$
{\mathcal C}_\beta(g_\gamma)(t)={\beta\over \Gamma(\gamma)t^{\beta-1}}\int_0^t(t-s)^{\beta-1}s^{\gamma -1}ds= { \Gamma(\beta+1)\Gamma(\gamma)\over \Gamma(\beta+\gamma)}g_\gamma(t), \qquad t>0.
$$

\noindent (ii) Take $e_\lambda(t):=e^{-\lambda t}$ for $t>0$ and $\lambda \in \CC^+$. Then
$$
{\mathcal C}_1(e_\lambda)(t)= {1\over \lambda t}(1-e^{-\lambda t}), \qquad {\mathcal C}_2(e_\lambda)(t)= {2\over \lambda t}(e^{-\lambda t}-1+\lambda t), \qquad t>0.
$$
Since ${\mathcal C}_1^2(e_\lambda)(t)= \displaystyle{{1\over t\lambda}\int_0^{t}{1-e^{-\lambda s}\over s}ds}$ for $t>0$, we conclude that ${\mathcal C}_1^2(e_\lambda)\not = {\mathcal C}_2(e_\lambda)$ and then ${\mathcal C}_1^2\not = {\mathcal C}_2$.

\noindent (iii) More generally, take $f_{\lambda}(t) := E_{\beta,1} (\lambda t^{\beta} )$ the Mittag-Leffler function, for
$t>0$ and $\lambda \in \CC^+$. Then
$$
{\mathcal C}_{\beta}(f_\lambda)(t)= {1\over \lambda g_{\beta + 1}(t)}(1-f_{\lambda}(t)), \quad t>0.
$$
}
\end{example}

The relationship between these generalized Ces\'{a}ro operators and fractional evolution equations of order $\alpha$ can be also
 observed in \cite{Li-Pr-03}.

The next Lemma shows a key commutativity property.

\begin{lemma}\label{llave}
Take $\alpha \ge 0 $ and $\beta>0$. Then $D^\alpha_+\circ \mathcal{C}_\beta=\mathcal{C}_\beta\circ D^\alpha_+$, i.e.,
$$
D^\alpha_+(\mathcal{C}_\beta(f))=\mathcal{C}_\beta( D^\alpha_+(f)),\qquad f\in {\mathcal S}_+,
$$
where $D^\alpha_+(t)=\displaystyle{1\over \Gamma(\alpha+1)}t^\alpha W^\alpha_+ f(t)$ for $f\in {\mathcal S}_+$.
\end{lemma}

\bgproof
By the equality (\ref{escala}), we have that
\begin{eqnarray*}
\mathcal{C}_\beta( D^\alpha_+(f))(t)&=&\beta\int_0^1(1-r)^{\beta-1}(tr)^\alpha W^\alpha_+f(tr)dr \cr
&=&t^\alpha W^\alpha_+\left(\beta \int_0^1(1-r)^{\beta-1} f(r)dr\right)=D^\alpha_+(\mathcal{C}_\beta(f))(t)
\end{eqnarray*}
for $f\in {\mathcal S}_+$ and we conclude the proof.
\edproof

The first main result in this section is the following theorem.

\begin{theorem}\label{lemma1.1}
The operator $\mathcal{C}_\beta$ is a bounded operator on $\mathcal{T}_p^{(\alpha)}(t^\alpha)$ and
$$
||\mathcal{C}_\beta||=\frac{\Gamma(\beta+1)\Gamma(1-1/p)}{\Gamma(\beta+1-1/p)},
$$
for $\alpha \ge 0$, $p>1$ and $\beta>0$. If $f\in \mathcal{T}_p^{(\alpha)}(t^\alpha)$, then
\begin{equation}\label{integral}
\mathcal{C}_\beta f(t)=\displaystyle\beta\int_0^\infty (1-e^{-r})^{\beta-1}e^{-r(1-1/p)}T_{r,p}f(t)dr, \quad t\geq 0,
\end{equation}
where the semigroup $(T_{r, p})_{t\ge 0}$ is defined in Theorem \ref{semigroup}.
\end{theorem}

\bgproof
We apply the change of variable $s= t e^{-r}$ to get that
$$
\mathcal{C}_\beta f(t) := \frac{\beta}{t^\beta}\int_0^t (t-s)^{\beta-1}f(s)ds = \beta\int_0^\infty (1-e^{-r})^{\beta-1}e^{-r}f(te^{-r})dr,
$$
and the equality (\ref{integral}) is proved. Observe that by this equality, $\mathcal{C}_\beta$ is well defined  and is a bounded operator on $\mathcal{T}_p^{(\alpha)}(t^\alpha)$ for $p>1:$ taken $f\in \mathcal{T}_p^{(\alpha)}(t^\alpha)$, then
\begin{eqnarray*}
||\mathcal{C}_\beta f||_{\alpha,p}&\leq& \beta\int_0^\infty (1-e^{-r})^{\beta-1}e^{-r(1-1/p)}||T_rf||_{\alpha,p}dr\\
&=&\beta||f||_{\alpha,p}\int_0^\infty (1-e^{-r})^{\beta-1}e^{-r(1-1/p)}dr=||f||_{\alpha,p}\frac{\Gamma(\beta+1)\Gamma(1-1/p)}{\Gamma(\beta+1-1/p)}.
\end{eqnarray*}
To check the exact value of $||\mathcal{C}_\beta ||$, note that by the Lemma \ref{llave} and the boudedness of $\mathcal{C}_\beta$ on $L^p(\RR^+)$  (see the Introduction) we get that
\begin{eqnarray*}
||\mathcal{C}_\beta f||_{\alpha,p}&=&||D^\alpha_+\circ \mathcal{C}_\beta f||_{p}= || \mathcal{C}_\beta\circ D^\alpha_+  f||_{p}\\&=& \frac{\Gamma(\beta+1)\Gamma(1-1/p)}{\Gamma(\beta+1-1/p)} || D^\alpha_+  f||_{p}= \frac{\Gamma(\beta+1)\Gamma(1-1/p)}{\Gamma(\beta+1-1/p)}||   f||_{\alpha,p},
\end{eqnarray*}
where we have applied Proposition \ref{juanjo} (iii).
\edproof

\begin{remark} {\rm (i) Recall that the Beta function, also called the Euler
integral of the first kind, is defined by:
$$
B(x,y)= \int_0^1 t^{x-1} (1-t)^{y-1} dt, \quad x>0, \quad y>0,
$$
and satisfies the property $ B(x,y)= \frac{\Gamma(x) \Gamma(y)}{\Gamma(x+y)}.$ Hence, the obtained value for the norm of $\mathcal{C}_\beta$ can be rewritten as
$$
||\mathcal{C}_\beta|| = \beta B(\beta, 1-1/p), \quad \beta >0, \quad p>1.
$$

(ii) In the case $p=1$ we remark that $\mathcal{C}_\beta$ does not take $\mathcal{T}_1^{(\alpha)}(t^\alpha)$ in $\mathcal{T}_1^{(\alpha)}(t^\alpha)$. In fact, from Lemma \ref{lemma2.3} it follows that, for $\beta>0$, $h_\beta(t):=(1+t)^{-(\beta+1)}$ belongs to $\mathcal{T}_1^{(\alpha)}(t^\alpha)$. By \cite[Formula 2, p.173]{Sa-Ki-Ma-93} and \cite[p. 38]{MOS}, we have
\begin{eqnarray*}
\mathcal{C}_\beta h_\beta(t)&=&\frac{\beta}{t^\beta}\int_0^t{(t-s)^{\beta-1}\over (1+s)^{\beta+1}}ds={}_2F_1(1,\beta+1;\beta+1;-t)=(1+t)^{-1},
\end{eqnarray*}
where $_2F_1$ denotes the Gaussian hypergeometric function,
$$
_2F_1(a,b;c;z):={\Gamma(c)\over \Gamma(b)\Gamma(a)}\sum_{n=0}^\infty{\Gamma(a+n)\Gamma(b+n)\over \Gamma(c+n)}{z^n\over n!}.
$$
Since $\mathcal{C}_\beta h_\beta$ does not belong to $L^1(\mathbb{R}^+)$ and $\mathcal{T}_1^{(\alpha)}(t^{\alpha})\hookrightarrow L^1(\mathbb{R}^+)$ (see Proposition \ref{juanjo} (i)), we obtain $\mathcal{C}_\beta h_\beta\not \in \mathcal{T}_1^{(\alpha)}(t^\alpha)$.

(iii) Take $\beta=1$ and  $f\in \mathcal{T}_p^{(\alpha)}(t^\alpha).$ Then
\begin{equation}\label{eq3.2}
\mathcal{C}_1 f(t)=\int_0^\infty e^{-r(1-1/p)}T_{r,p}f(t)dr=R(\lambda_p,\Lambda)f(t),\quad \lambda_p=1-1/p>0.
\end{equation}
and by the spectral theorem for resolvent operators (see for example \cite[Theorem IV.1.13]{En-Na-00}) we get that
\begin{equation}\label{eq3.3}
\sigma(\mathcal{C}_1)=\left\{w\in\mathbb{C} \ : \ \left|w-\frac{p}{2(p-1)}\right|=\frac{p}{2(p-1)}\right\},
\end{equation}
see \cite[Theorem 2]{Mo99} and similar results in \cite[Theorem 3.1]{Ar-Si-11},  and \cite[Corollary 2.2]{Ar-Pa-02}.
Here, $R(\cdot,\Lambda)$ denotes the resolvent operator of $\Lambda.$

Note that in case $\beta =2$ we obtain
\begin{equation*}
\mathcal{C}_2 f(t)= 2 \int_0^\infty e^{-r(1-1/p)}(1-e^{-r})T_{r,p}f(t)dr=2 R(\lambda_p,\Lambda)f(t) - 2 R(\lambda_p +1 , \Lambda )f(t),
\end{equation*}
and, more generally, for $\beta =n+1$,
\begin{equation}
\mathcal{C}_{n+1} f(t)= (n+1)\sum_{k=0}^n  \binom{n}{k} (-1)^k R(\lambda_p + k, \Lambda )f(t), \quad  n \in \mathbb{Z}_+.
\end{equation}
}
\end{remark}

In the next result, we are able to describe $\sigma(\mathcal{C}_\beta)$ for $\beta>0$.

\begin{theorem}\label{spec}
Let $1<p<\infty,$ and $\mathcal{C}_\beta :\mathcal{T}_p^{(\alpha)}(t^\alpha) \to \mathcal{T}_p^{(\alpha)}(t^\alpha)$ the generalized Ces\`{a}ro operator. Then

$$
\sigma(\mathcal{C}_\beta)= \beta \overline{B(\beta, 1-1/p + i \mathbb{R})}:= \Gamma(\beta+1)\overline{\left\{{\Gamma(1-{1\over p}+it)\over \Gamma(\beta+1-{1\over p}+it)} \ : \ t\in \mathbb{R}\right\}}.
$$
\end{theorem}

\bgproof
As $(T_{t,p})_{t\in \RR}$ is an uniformly bounded $C_0$-group (Theorem \ref{semigroup}) whose infinitesimal generator is $(\Lambda, D(\Lambda))$ and $\mathcal{C}_\beta= \widehat{f_{\beta, p}}(\Lambda)$, i.e,
$$
\mathcal{C}_\beta f=\displaystyle\beta\int_0^\infty (1-e^{-r})^{\beta-1}e^{-r(1-1/p)}T_{r,p}fdr= \int_{-\infty}^\infty f_{\beta,p}(r)T_{r,p}fdr,
$$
where $f_{\beta,p}(r)= \chi_{[0,\infty)}(r)\beta (1-e^{-r})^{\beta-1}e^{-r(1-1/p)}$ for $r\in \RR$, see  Theorem \ref{lemma1.1}. By  \cite[Theorem 3.1]{[Se]}, we obtain
$$
\sigma(\mathcal{C}_\beta)=\overline{\widehat{f_{\beta, p}}(\sigma(i\Lambda))}
$$
where $\widehat{f_{\beta, p}}$ is the Fourier transform of the function $f_{\beta, p}$.
As $\sigma(i\Lambda)=\RR$ (see Proposition \ref{spectrum} (ii)) and $\widehat{f_{\beta, p}}(t)={\mathcal L}(f_{\beta, p})(it)$ we use that
$$
{\mathcal L}(f_{\beta, p})(z)=\beta\int_0^\infty e^{-zr}(1-e^{-r})^{\beta-1}e^{-r(1-1/p)}dr= {\Gamma(\beta+1)\Gamma(1-{1\over p}+z)\over \Gamma(\beta+1-{1\over p}+z)}, \qquad z\in \overline{\CC^+}.
$$
to conclude the result.
\edproof

\begin{remark}
{\rm  In the case that $n\in\NN$, we obtain that
\begin{eqnarray*}
\sigma(\mathcal{C}_n)
&=&\left\{{n!p^n\over ((n+it)p-1)\dots((1+it)p-1)} \ : \ t\in \mathbb{R}\right\}\cup\{0\},
\end{eqnarray*}
and for $n=1$
\begin{eqnarray*}
\sigma(\mathcal{C}_1)
&=&\left\{{p\over (1+it)p-1} \ : \ t\in \mathbb{R}\right\}\cup\{0\}=\left\{w\in\mathbb{C} \ : \ \left|w-\frac{p}{2(p-1)}\right|=\frac{p}{2(p-1)}\right\}.
\end{eqnarray*}
}
\end{remark}

Now we consider the generalized dual Ces\`{a}ro operator $\mathcal{C}_\beta^*$  on $\mathcal{T}_p^{(\alpha)}(t^\alpha)$ defined by
$$
\mathcal{C}_\beta^* f(t) := {\beta}\int_t^\infty {(s-t)^{\beta-1}\over s^\beta} f(s)ds = \beta\int_1^\infty{(r-1)^{\beta-1}\over r^\beta}f(tr)dr, \quad \,t> 0.
$$

For $0<\gamma<1$, functions  $g_\gamma$ are eigenfunctions of ${\mathcal C}_\beta^*$ with eigenvalue ${\Gamma(\beta+1) \Gamma(1-\gamma)\over \Gamma(\beta-\gamma+1)}$:
$$
{\mathcal C}_\beta^*(g_\gamma)(t)={\beta\over \Gamma(\gamma)}\int_t^\infty{(s-t)^{\beta-1}s^{\gamma -1}\over s^{\beta}}ds={\Gamma(\beta+1) \Gamma(1-\gamma)\over \Gamma(\beta-\gamma+1)}g_\gamma(t),
$$
for $t>0$.

Using similar ideas to Lemma \ref{llave}, we obtain
\begin{equation}\label{conmu}
D^\alpha_+\circ \mathcal{C}_\beta^*(f) = \mathcal{C}_\beta^*\circ D^\alpha_+(f),\qquad f\in {\mathcal S}_+
\end{equation}
where $D^\alpha_+ f(t)=\displaystyle{1\over \Gamma(\alpha+1)}t^\alpha W^\alpha_+f(t)$ for $f\in {\mathcal S}_+$ and $t\ge 0$. Hence we are ready to prove our main dual result.

\begin{theorem}\label{theorem3.6}
The operator $\mathcal{C}_\beta^*$ is a bounded operator on $\mathcal{T}_p^{(\alpha)}(t^\alpha)$ and
$$
||\mathcal{C}_\beta^*||=\frac{\Gamma(\beta+1)\Gamma(1/p)}{\Gamma(\beta+1/p)},
$$
for $\alpha \ge 0$, $p > 1$ and $\beta>0$. The dual operator of $\mathcal{C}_\beta$ on $\mathcal{T}_p^{(\alpha)}(t^\alpha)$ is $\mathcal{C}_\beta^*$  on $\mathcal{T}_{p'}^{(\alpha)}(t^\alpha)$, i.e.
$$
\langle \mathcal{C}_\beta f,g \rangle_\alpha=\langle f,\mathcal{C}_\beta^*g \rangle_\alpha, \qquad f \in \mathcal{T}_p^{(\alpha)}(t^\alpha), \quad  g \in\mathcal{T}_{p'}^{(\alpha)}(t^\alpha),
$$
where $\langle \quad,\quad \rangle_\alpha$ is given in Proposition \ref{juanjo} (iv) and $\frac{1}{p}+\frac{1}{p'}=1$.

If $f\in \mathcal{T}_p^{(\alpha)}(t^\alpha)$, then
\begin{eqnarray*}\label{integral2}
\mathcal{C}_\beta^* f(t)&=& \displaystyle\beta\int_{-\infty}^0 (e^{-r}-1)^{\beta-1}e^{-r(1-1/p-\beta)}T_{r,p}f(t)dr, \quad t\geq 0,\end{eqnarray*}
where the $C_0$-group $(T_{r, p})_{t\in \RR}$ is defined in Theorem \ref{semigroup}.
\end{theorem}

\bgproof
We change the variable $s = t e^r$ to obtain
\begin{eqnarray*}
\mathcal{C}_\beta^* f(t)&=&{\beta}\int_t^\infty {(s-t)^{\beta-1}\over s^\beta} f(s)ds={\beta}\int_0^\infty (e^r-1)^{\beta-1} e^{r(1-\beta)} f(te^r)dr\cr
&=& {\beta}\int_0^\infty (e^r-1)^{\beta-1} e^{r(1-\beta-{1\over p})}e^{r\over p} f(te^r)dr\cr
&=& \beta\int_{-\infty}^0 (e^{-r}-1)^{\beta-1}e^{-r(1-1/p-\beta)}T_{r,p}f(t)dr
\end{eqnarray*}
for $ f\in \mathcal{T}_p^{(\alpha)}(t^\alpha)$, $p\ge 1$ and $t\ge 0$. To check $||\mathcal{C}_\beta^* ||$, note that by the Lemma \ref{llave} and the boudedness of $\mathcal{C}_\beta^*$ on $L^p(\RR^+)$  (see the Introduction) we get that
\begin{eqnarray*}
||\mathcal{C}_\beta^* f||_{\alpha,p}&=&||D^\alpha_+\circ \mathcal{C}_\beta^* f||_{p}= || \mathcal{C}_\beta^*\circ D^\alpha_+  f||_{p}\\&=& \frac{\Gamma(\beta+1)\Gamma(1/p)}{\Gamma(\beta+1/p)} || D^\alpha_+  f||_{p}= \frac{\Gamma(\beta+1)\Gamma(1/p)}{\Gamma(\beta+1/p)}||   f||_{\alpha,p},
\end{eqnarray*}
where we have applied the Proposition \ref{juanjo} (iii) and the equality (\ref{conmu}).

To conclude the proof, take $f \in \mathcal{T}_p^{(\alpha)}(t^\alpha),$ and  $g \in\mathcal{T}_{p'}^{(\alpha)}(t^\alpha)$. Then
\begin{eqnarray*}
\langle \mathcal{C}_\beta f,g \rangle_\alpha&=&\langle D^\alpha_+(\mathcal{C}_\beta f),D^\alpha_+ g \rangle_0=\langle \mathcal{C}_\beta(D^\alpha_+ f),D^\alpha_+ g \rangle_0=\langle  D^\alpha_+ f,\mathcal{C}_\beta^\ast(D^\alpha_+ g )\rangle_0\cr
&=& \langle  D^\alpha_+ f,D^\alpha_+(\mathcal{C}_\beta^\ast g )\rangle_0= \langle f,\mathcal{C}_\beta^\ast g )\rangle_\alpha,
\end{eqnarray*}
where we have used  Lemma \ref{llave},  $(\mathcal{C}_\beta)^\ast= \mathcal{C}_\beta^\ast$ and equality (\ref{conmu}).
\edproof

\begin{remark} {\rm Take $\beta=1$ and $f\in
\mathcal{T}_p^{(\alpha)}(t^\alpha).$ Then
\begin{equation*}
\mathcal{C}_1^* f(t)=\int_{-\infty}^0 e^{-{r\over p}}T_{-r,p}f(t)drds=R(1/p,-\Lambda)f(t), \qquad t\ge 0.
\end{equation*}
and by the spectral theorem for the resolvent operator, see \cite[Theorem IV.1.13]{En-Na-00}, we obtain
\begin{equation*}
\sigma(\mathcal{C}_1^*)=\left\{w\in\mathbb{C} \ : \ \left|w-\frac{p}{2}\right|=\frac{p}{2}\right\}.
\end{equation*}
This  gives a proof of a conjecture posed by F. M\'{o}ricz in
\cite[Section 2]{Mo99}. See a similar result in \cite[Theorem
3.2]{Ar-Si-11}.}

\end{remark}
 In the following theorem we describe $\sigma(
\mathcal{C}_\beta^*)$  for $\beta>0$.

\begin{theorem}\label{spectrr}
Let $\beta>0$, $1\le p<\infty,$ and $\mathcal{C}_\beta^\ast :\mathcal{T}_p^{(\alpha)}(t^\alpha) \to \mathcal{T}_p^{(\alpha)}(t^\alpha)$ the generalized dual Ces\`{a}ro operator. Then

$$
\sigma(\mathcal{C}_\beta^\ast)= \beta \overline{B(\beta, 1/p + i\mathbb{R})}:= \Gamma(\beta+1)\overline{\left\{{\Gamma({1\over p}+it)\over \Gamma(\beta+{1\over p}+it)} \ : \ t\in \mathbb{R}\right\}}.
$$
\end{theorem}

\bgproof
We remind that $(T_{t,p})_{t\in \RR}$ is a uniformly bounded $C_0$-group  whose infinitesimal generator is $(\Lambda, D(\Lambda))$ and
$\mathcal{C}_\beta^\ast= \widehat{g_{\beta, p}}(\Lambda)$, i.e,
$$
\mathcal{C}_\beta^\ast f=\displaystyle\beta\int_{-\infty}^0 (e^{-r}-1)^{\beta-1}e^{-r(1-1/p-\beta)}T_{r,p}fdr= \int_{-\infty}^\infty g_{\beta,p}(r)T_{r,p}fdr,
$$
where $g_{\beta,p}(r)= \chi_{(-\infty,0]}(r)\beta (e^{-r}-1)^{\beta-1}e^{-r(1-1/p-\beta)}$ for $r\in \RR$, see  Theorem \ref{theorem3.6}. By  \cite[Theorem 3.1]{[Se]}, we obtain
$$
\sigma(\mathcal{C}_\beta^\ast)=\overline{\widehat{g_{\beta, p}}(\sigma(i\Lambda))}
$$
where $\widehat{g_{\beta, p}}$ is the Fourier transform of the function $g_{\beta, p}$.
As $\sigma(i\Lambda)=\RR$, and  $\widehat{{g_{\beta, p}}}(t)={\mathcal L}(\widetilde{g_{\beta, p}})(-it)$ (where $\widetilde{g_{\beta, p}}(\xi):= g_{\beta, p}(-\xi)$ for $\xi \ge 0$) we use that
$$
{\mathcal L}(\widetilde{g_{\beta, p}})(z)=\beta\int_0^\infty e^{-zr}(e^{r}-1)^{\beta-1}e^{r(1-1/p-\beta)}dr= {\Gamma(\beta+1)\Gamma({1\over p}+z)\over \Gamma(\beta+{1\over p}+z)}, \qquad z\in \overline{\CC^+},
$$
to conclude the result.
\edproof

\begin{remark}
{\rm  In the case that $n\in \NN$, we obtain that
\begin{eqnarray*}
\sigma(\mathcal{C}_n^\ast)
&=&\left\{{n!p^n\over ((n-1)p+1+it)\dots(p+1+it)(1+it)} \ : \ t\in \mathbb{R}\right\}\cup\{0\},
\end{eqnarray*}
and for $n=1$
\begin{eqnarray*}
\sigma(\mathcal{C}_1^*)
&=&\left\{{p\over 1+it} \ : \ t\in \mathbb{R}\right\}\cup\{0\}=\left\{w\in\mathbb{C} \ : \ \left|w-\frac{p}{2}\right|=\frac{p}{2}\right\}.
\end{eqnarray*}
}
\end{remark}

\begin{remark}
{\rm In the case that $p=2$ we have
$\sigma(\mathcal{C}_{\beta})= \sigma(\mathcal{C}_{\beta}^*)$ for all $\beta >0$. Note that in case $p\neq 2$ the spectrum of $\mathcal{C}_{\beta}$ and $\mathcal{C}_{\beta}^*$ are dual in the sense that $\sigma(\mathcal{C}_{\beta}),$ with $\mathcal{C}_{\beta}$ defined on $\mathcal{T}_p^{(\alpha)}(t^\alpha),$ is identical to $\sigma(\mathcal{C}_{\beta}^*),$ with $\mathcal{C}_{\beta}^*$ defined on
$\mathcal{T}_{p'}^{(\alpha)}(t^\alpha)$, and where $\frac{1}{p} + \frac{1}{p'}=1.$
}
\end{remark}

To finish this section we prove the remarkable fact that
$\mathcal{C}_\alpha$ and $\mathcal{C}_\beta^\ast$ commute on
$L^p(\R^+)$ (and then on $\mathcal{T}_p^{(\alpha)}(t^\alpha)$). We
also give explicitly the value of
$\mathcal{C}_\alpha\mathcal{C}_\beta^\ast$ in terms of the  the
Gaussian hypergeometric function $_2F_1$. This theorem includes
\cite[Lemma 2]{Mo99} for $\alpha=\beta=1$.

\begin{theorem} \label{conmutacion} Let $\mathcal{C}_\alpha$ and $\mathcal{C}_\beta^\ast$ the generalized Ces\'{a}ro operators on $L^p(\R^+)$ for $p>1$.
Then $\mathcal{C}_\alpha\mathcal{C}_\beta^\ast=\mathcal{C}_\beta^\ast\mathcal{C}_\alpha$ for $\alpha, \beta >0$ and
\begin{eqnarray*}
(\mathcal{C}_\alpha\mathcal{C}_\beta^\ast)f(t)&=&\alpha\int^t_0f(r){1\over t-r}\left({t-r\over t}\right)^{\alpha+\beta}
\,_2F_1(\alpha+\beta, \beta; \beta+1; {r\over t})dr\cr\qquad &\quad&\qquad + \,\, \beta\int_t^\infty f(r){1\over r-t}\left({r-t\over t}\right)^{\alpha+\beta}
\,_2F_1(\alpha+\beta,  \alpha; \alpha+1; {t\over r})dr,
\end{eqnarray*}
 in particular
\begin{eqnarray*}
(\mathcal{C}_1\mathcal{C}_\beta^\ast )f(t)&=&\mathcal{C}_1 f(t)+\beta\int_t^\infty f(r){(r-t)^\beta\over r^{\beta+1}}
\,_2F_1(\beta+1, 1; 2; {r\over t})dr,\\
(\mathcal{C}_\alpha\mathcal{C}_1^\ast )f(t)&=&\alpha\int^t_0f(r){(t-r)^{\alpha}\over t^{\alpha+1}}
\,_2F_1(\alpha+1, 1; 2; {r\over t})dr+ \mathcal{C}_1^\ast f(t),\\
(\mathcal{C}_1\mathcal{C}_1^\ast )f&=&\mathcal{C}_1f+ \mathcal{C}_1^\ast f=(\mathcal{C}^*_1\mathcal{C}_1)f,
\end{eqnarray*}
for $f\in L^p(\R^+)$ and $t$ almost everywhere on $\R^{+}$.
\end{theorem}

\bgproof
By the integral representations (\ref{integral}) and (\ref{integral2}), and since $T_{t,p}$ commutes  with  $T_{r,p}$ for any $t, r\in \R$, we conclude that  $\mathcal{C}_\alpha\mathcal{C}_\beta^\ast=\mathcal{C}_\beta^\ast\mathcal{C}_\alpha$ for $\alpha, \beta >0$. Take $f \in L^p(\R^+)$ and we apply the Fubini theorem to get that
\begin{eqnarray*}
\mathcal{C}_\beta^\ast \mathcal{C}_\alpha f(t)&=&\beta\alpha\int_t^\infty {(x-t)^{\beta-1}\over x^{\beta+\alpha}}\int_0^x(x-r)^{\alpha-1}f(r) drdx\cr
&=& \beta\alpha\int_0^\infty f(r)\int_{\max\{t,r\}}^\infty {(x-t)^{\beta-1}(x-r)^{\alpha-1}\over x^{\beta+\alpha}}dxdr
\end{eqnarray*}
for $t$ almost everywhere on $\R^{+}$. For $0<r<t$, this equality
$$
\int_{t}^\infty {(x-t)^{\beta-1}(x-r)^{\alpha-1}\over x^{\beta+\alpha}}dx={1\over \beta(t-r)}\left({t-r\over t}\right)^{\alpha+\beta}\,_2F_1(\alpha+\beta, \beta; \beta+1; {r\over t})
$$
holds, see for example \cite[p. 314, 3197(1)]{[GR]}.

Now take $\alpha=1$. Since
$$
(1-z)^{a}\,_2F_1(a, b;c; z)= \,_2F_1(a, c-b;c; {z\over z-1})
$$
(see for example \cite[p.47]{MOS}), we get that
$$
{1\over t-r}\left({t-r\over t}\right)^{1+\beta}
\,_2F_1(1+\beta, \beta; \beta+1; {r\over t})= {1\over t-r}\,_2F_1(1+\beta, 1;1+\beta; {-r\over t-r})={1\over t}
$$
where we apply that $\,_2F_1(-a, b;b; -z)=(1+z)^a$, (\cite[p. 38]{MOS}). Similarly we prove the case $\beta=1$.
\edproof

\section{Composition groups on Sobolev spaces defined on $\mathbb{R}.$}

\setcounter{theorem}{0}
\setcounter{equation}{0}

In  this section we introduce the subspaces $\mathcal{T}_p^{(\alpha)}(|t|^\alpha)$ which are contained in $L^p(\RR)$, similarly to   $\mathcal{T}_p^{(\alpha)}(t^\alpha)$ are in $L^p(\RR^+)$. Let ${\mathcal S}$ be the Schwartz class on $\mathbb{R}$ and we set
$$
W_{-}^{-\alpha}f(x)={1\over
\Gamma(\alpha)}\int_{-\infty}^{x}(x-t)^{\alpha-1}f(t)dt,
$$
$$
W_{-}^{\alpha}f(x)={1\over\Gamma(n-\alpha)}{d^n \over
dx^n}\int_{-\infty}^{x}(x-t)^{n-\alpha-1}f(t)dt,
$$
and $W_{-}^0f=f$, for $x\in\R$ and a natural number $n>\alpha$. Putting $\tilde{f}(x)=f(-x)$, it is readily seen that $W_{+}^{\alpha}f(x) = W_{-}^{\alpha}\tilde{f}(-x)$ for all $\alpha\in\R$, $f\in {\mathcal S}$ and $x\in\R$. Equalities $W_{-}^{\alpha+\beta}=W_{-}^\alpha W_{-}^\beta$ and $W_-^nf=f^{(n)}$ hold for each natural number $n$ and $\alpha,\beta\in\R$.

For $f\in {\mathcal S}$, put
$$
W^{\alpha}_0f(t):= \left\{\begin{array}{ll}   W^\alpha_-f(t),              & t<0,\\
                                              e^{i\pi\alpha}W^\alpha_+f(t),& t>0.\\
                          \end{array} \right.
$$
For $\lambda >0$, we have that $W_0^\alpha (f_\lambda)= \lambda^\alpha (W_0^\alpha f)_\lambda$, where $f_\lambda(t)= f(\lambda t)$ for $t\in \R$.

\begin{definition}
Let $1\le p<\infty$. The Banach space $\mathcal{T}_p^{(\alpha)}(|t|^\alpha)$  is defined  as the completion of the Schwartz class on $\mathbb{R}$ in the norm
\begin{equation*}
|||f|||_{\alpha, p}:={1\over \Gamma( \alpha+1)}\left(\int_{-\infty}^\infty \left(|W^\alpha_0 f(t)\vert \,\vert t|^{\alpha }\right)^pdt\right)^{\frac{1}{p}}.
\end{equation*}
\end{definition}

Similar properties to $\mathcal{T}_p^{(\alpha)}(t^\alpha)$ hold in $\mathcal{T}_p^{(\alpha)}(|t|^\alpha)$. The proof of next proposition is similar to the proof of Proposition \ref{juanjo} and we skip it.

\begin{proposition}\label{juanjo2}
Take  $p\geq 1$ and $\beta>\alpha>0$. Then
\begin{itemize}
\item[(i)] $\mathcal{T}_p^{(\beta)}(\vert t\vert^\beta)\hookrightarrow \mathcal{T}_p^{(\alpha)}(\vert t\vert^\alpha)\hookrightarrow  L^p(\mathbb{R}) $.

\item[(ii)] The operator $D^\alpha_0: \mathcal{T}_p^{(\alpha)}(\vert t\vert^\alpha)\to L^p(\RR)$  defined by
$$
f\mapsto D^\alpha_0 f(t):={1\over \Gamma(\alpha+1)}\vert t\vert^\alpha W^\alpha_0f(t), \qquad t\in \R, \quad f \in \mathcal{T}_p^{(\alpha)}(\vert t\vert^\alpha),
$$
is an isometry.

\item[(iii)] If $p>1$ and $p'$ satisfies $\frac{1}{p}+\frac{1}{p'}=1$, then the dual of $\mathcal{T}_p^{(\alpha)}(\vert t\vert^\alpha)$ is $\mathcal{T}_{p'}^{(\alpha)}(\vert t\vert^\alpha)$, where the duality is given by
$$
\langle f,g \rangle_{\alpha} ={1\over \Gamma(\alpha+1)^2}\int_{-\infty}^\infty W^\alpha_0 f(t)W^\alpha_0 g(t)\vert t\vert^{2\alpha}dt,
$$
for $f\in \mathcal{T}_p^{(\alpha)}(\vert t\vert^\alpha)$, $g\in \mathcal{T}_{p'}^{(\alpha)}(\vert t\vert^\alpha)$.
\end{itemize}
\end{proposition}

For $p=1$, the subspace $\mathcal{T}_1^{(\alpha)}(|t|^\alpha)$ was introduced in \cite[Definition 1.9]{Ga-Mi-06}. In fact $\mathcal{T}_1^{(\alpha)}(|t|^\alpha)$  is a subalgebra of $L^1(\R)$ for the convolution product
\begin{equation}\label{convos2}
f\ast g(t)=\int_{-\infty}^{\infty} f(t-s)g(s)ds, \qquad t\in \R,\quad  f, g\in \mathcal{T}_1^{(\alpha)}(|t|^\alpha),
\end{equation}
see \cite[Theorem 1.8]{Ga-Mi-06} and also \cite[Theorem 2]{Mi-07} for some more details.

\begin{theorem}\label{module}
Let $1< p<\infty$. The Banach space $\mathcal{T}_p^{(\alpha)}(|t|^\alpha)$  is a module for the algebra $\mathcal{T}_1^{(\alpha)}(|t|^\alpha)$ and
$$
|||f\ast g|||_{\alpha, p}\le C_{\alpha, p} |||f|||_{\alpha, p}||| g|||_{\alpha, 1}, \qquad f\in \mathcal{T}_p^{(\alpha)}(|t|^\alpha), \quad g \in \mathcal{T}_1^{(\alpha)}(|t|^\alpha).
$$
\end{theorem}

\bgproof
Take $f,g\in \S$. We  write  $f_+:=f\chi_{[0,\infty)}$ and $f_{-}:=f\chi_{(-\infty,0]}$. By considering the  decomposition $f\ast g = (f_+ \ast g_+) + (f_+ \ast g_-) + (f_-\ast g_+) + (f_-\ast g_-)$ on $\R$, and we apply \cite[Lemma 1.6]{Ga-Mi-06} and the fact that $f_-\ast g_-=0$ on $(0,\infty)$ to obtain that
$$
W^{\alpha}_+(f\ast g)_+(t) =W^{\alpha}_+(f_+\ast
g_+)(t)+(W^{\alpha}_+f_+\ast g_-)(t) +(W^{\alpha}_+g_+\ast
f_-)(t),  \qquad t>0.
$$
Now, first,
$$
\Vert f_+\ast g_+\Vert_{\alpha, p}\le C_{\alpha, p}\Vert f_+\Vert_{\alpha, p}\Vert  g_+\Vert_{\alpha, 1}\le C_{\alpha, p}||| f|||_{\alpha, p}||| g|||_{\alpha, 1}
$$
by Proposition \ref{juanjo} (ii).

On the other hand, $\mathcal{T}_1^{(\alpha)}(t^\alpha) \subset L^1(\R^+)$, and we apply the Minkowski inequality to get that
\begin{eqnarray*}
&\,&\left(\int_0^\infty\vert W^{\alpha}_+f_+\ast g_-(t)\vert^p t^{\alpha p} dt\right)^{1\over p}
\le\left(\int_0^\infty\left(\int_0^\infty \vert W^{\alpha}_+f_+(s+t)\vert\vert
g_-(s)\vert ds\right)^p t^{\alpha p} dt\right)^{1\over p}\\&=&\int_0^\infty\vert
g_-(s)\vert \left(\int_0^\infty \vert W^{\alpha}_+f_+(t+s)\vert^p t^{\alpha p} dt\right)^{1\over p}ds
\le
\int_0^\infty\vert
g_-(s)\vert \left(\int_s^\infty \vert W^{\alpha}_+f_+(u)\vert^p u^{\alpha p} du\right)^{1\over p}ds\\
& \le & \Gamma(\alpha+1)|||g|||_{0,1} || f_+||_{\alpha, p}\le \Gamma(\alpha+1)|||g|||_{\alpha,1} ||| f|||_{\alpha, p}.
\end{eqnarray*}

As $\mathcal{T}_p^{(\alpha)}(t^\alpha)\subset L^p(\R^+)$ for $p>1$, and we apply again the Minkowski inequality to obtain that
\begin{eqnarray*}
&\,&\left(\int_0^\infty\vert (W^{\alpha}_+g_+\ast
f_-)(t)(t)\vert^p t^{\alpha p} dt\right)^{1\over p}
\le\left(\int_0^\infty\left(\int_t^\infty \vert W^{\alpha}_+g_+(s)\vert\vert
f_-(t-s)\vert ds\right)^p t^{\alpha p} dt\right)^{1\over p}\\&=&\int_0^\infty\vert
W^{\alpha}_+g_+(s)\vert \left(\int_0^s \vert f_-(t-s)\vert^p t^{\alpha p} dt\right)^{1\over p}ds
\le
||| f|||_{0,p}\int_0^\infty\vert
W^{\alpha}_+g_+(s)\vert s^\alpha ds\\
& \le & \Gamma(\alpha+1)|||f|||_{\alpha,p}\, || g_+||_{\alpha, 1}\le \Gamma(\alpha+1)|||f|||_{\alpha,p}\, ||| g|||_{\alpha, 1}.
\end{eqnarray*}

Combining these estimates obtained, we get
$$
{1\over \Gamma(\alpha+1)}\left(\int_0^\infty\vert W^{\alpha}_+(f\ast g)(t)\vert^p\ t^{\alpha p} dt\right)^{1\over p}
\le C||| f|||_{\alpha,p}\,|||
g|||_{\alpha,1}.
$$
Finally, because $W^{\alpha}_-(f\ast g)(t)=W^{\alpha}_+(\tilde{f}\ast\tilde{g})(-t)$ if $t<0$ using the inclusion
${\mathcal T}_p^{(\alpha)}(t^\alpha)\subset L^p(\R^+)$ as above for $p\ge 1$, we have that
$$
{1\over \Gamma(\alpha+1)}\left(\int_{-\infty}^0\vert W^{\alpha}_-(f\ast g)(t)\vert^p \,\vert t\vert^{\alpha p}\
dt\right)^{1\over p} \le C ||| f|||_{\alpha,p}\, |||
g|||_{\alpha,1}.
$$
The result follows.
\edproof

We remark that, as in the case of $\mathcal{T}_p^{(\alpha)}(t^\alpha)$, it is easy to verify that $(T_{t,p})_{t\in \mathbb{R}}$ is a $C_0$-group of isometries on $\mathcal{T}_p^{(\alpha)}(|t|^\alpha)$ as the next theorem shows. The proof runs parallel to the proofs of Theorem \ref{semigroup}, Proposition \ref{spectrum} and Proposition \ref{dualsemi} and hence we omit it.

\begin{theorem}\label{groups}
Let $1\le p$ and $\alpha\ge 0$. We define the family of operators  $(T_{t,p})_{t\in \RR}$ by
$$
T_{t,p} f(s):=e^{-\frac{t}{p}}f(e^{-t}s), \qquad f \in \mathcal{T}_p^{(\alpha)}(\vert t\vert^\alpha).
$$
\begin{itemize}
\item[(i)] Then $(T_{t,p})_{t\in \RR}$ is a $C_0$-group of isometries on $\mathcal{T}_p^{(\alpha)}(\vert t\vert^\alpha)$ whose infinitesimal generator $\Lambda$ is given by
$$
(\Lambda f)(s):=-sf'(s)-\frac{1}{p}f(s)
$$
with domain $D(\Lambda)= \mathcal{T}_p^{(\alpha+1)}(\vert t\vert^{\alpha+1})$.

\item[(ii)] $\sigma_p(\Lambda)=\emptyset$ and $\sigma(\Lambda)=i\mathbb{R}$ (here $\sigma_p$ denotes the point spectrum).

\item[(iii)] The semigroups $(T_{t,p})_{t\geq 0}$ and $(T_{-t,p'})_{t\geq 0}$ are dual operators of each other acting on $\mathcal{T}_p^{(\alpha)}(\vert t\vert^\alpha)$ and $\mathcal{T}_{p'}^{(\alpha)}(\vert t\vert^\alpha)$ with  ${1\over p}+ {1\over p'}=1$ for $p>1.$
\end{itemize}
\end{theorem}

\section{The generalized Ces\`{a}ro operators on $\mathbb{R}.$}

\setcounter{theorem}{0}
\setcounter{equation}{0}

For $\beta > 0$ we define the generalized Ces\`{a}ro operator by
$$
\mathcal{C}_\beta f(t) := \left\{ \begin{array}{ll} \displaystyle\frac{\beta}{|t|^\beta}\int_t^0 (s-t)^{\beta-1}f(s)ds, & t<0, \\
                                                    \\
                                                    f(0),                                                               & t=0, \\
                                                    \\
                                                    \displaystyle\frac{\beta}{t^\beta}\int_0^t (t-s)^{\beta-1}f(s)ds,   & t>0, \\
                                  \end{array} \right.
$$
for $f\in {\mathcal S}.$ We are interested in the extension of  $\mathcal{C}_\beta$ on $\mathcal{T}_p^{(\alpha)}(\vert t\vert^\alpha)$. Note that we may write
$$
\mathcal{C}_\beta f(t)= \beta\int_0^1(1-r)^{\beta-1}f(tr)dr, \quad \,\,t\in \R, \,\, f\in {\mathcal S}.
$$
We use this integral representation to proof the next lemma.

\begin{lemma}\label{llave2}
Take $\alpha \ge 0 $ and $\beta>0$. Then $ D^\alpha_0\circ \mathcal{C}_\beta=\mathcal{C}_\beta\circ D^\alpha_0,$ i.e.,
$$
D^\alpha_0(\mathcal{C}_\beta(f))=\mathcal{C}_\beta( D^\alpha_0(f)),\qquad f\in {\mathcal S},
$$
where $D^\alpha_0f(t)=\frac{1}{\Gamma(\alpha+1)}\vert  t\vert
^\alpha W^\alpha_0 f(t)$ for $f\in {\mathcal S}$.
\end{lemma}

\bgproof
Since for $\lambda >0$, we have that $W_0^\alpha (f_\lambda)= \lambda^\alpha (W_0^\alpha f)_\lambda$, where $f_\lambda(t)= f(\lambda t)$ for $t\in \R,$ the proof follows similarly to Lemma \ref{llave}.
\edproof

Similar results of $\mathcal{C}_\beta$ on $\mathcal{T}_p^{(\alpha)}(t^\alpha)$ hold for $\mathcal{C}_\beta$ on $\mathcal{T}_p^{(\alpha)}(\vert t\vert^\alpha)$. The proof of next result is analogous to the proof of Theorem \ref{lemma1.1} and Theorem \ref{spec}.

\begin{theorem}\label{lemma2.1}
Let $\alpha\geq 0$, $\beta > 0$, $1<p<\infty$ and the generalized Ces\`{a}ro operator $\mathcal{C}_\beta$ on $\mathcal{T}_p^{(\alpha)}(\vert t\vert^\alpha)$. Then
\begin{itemize}
\item[(i)] The operator $\mathcal{C}_\beta$ is bounded  on $\mathcal{T}_p^{(\alpha)}(\vert t\vert^\alpha)$ and
$$
||\mathcal{C}_\beta||=\frac{\Gamma(\beta+1)\Gamma(1-1/p)}{\Gamma(\beta+1-1/p)}.
$$

\item[(ii)] If $f\in \mathcal{T}_p^{(\alpha)}(\vert t\vert^\alpha)$, then
$$
\mathcal{C}_\beta f(t)=\displaystyle\beta\int_0^\infty (1-e^{-r})^{\beta-1}e^{-r(1-1/p)}T_{r,p}f(t)dr, \quad t\in \R ,
$$
where the $C_{0}$-group $(T_{r, p})_{ r\in\R}$ is defined in Theorem \ref{groups}.

\item[(iii)]
$$
\sigma(\mathcal{C}_\beta)=\Gamma(\beta+1)\overline{\left\{{\Gamma(1-{1\over p}+it)\over \Gamma(\beta+1-{1\over p}+it)} \ : \ t\in \R \right\}}.
$$
\end{itemize}
\end{theorem}

Now we consider the generalized dual Ces\`{a}ro operator $\mathcal{C}^*_\beta$ defined for $\beta > 0$ by
$$
\mathcal{C}^*_\beta f(t):= \left\{ \begin{array}{ll} \displaystyle\beta\int_{-\infty}^t \frac{(t-s)^{\beta-1}}{|s|^\beta}f(s)ds, & t<0, \\
                                                     \\
                                                     0,                                                                          & t=0, \\
                                                     \\
                                                     \displaystyle\beta\int_t^\infty \frac{(s-t)^{\beta-1}}{s^\beta}f(s)ds,      & t>0, \\
                                   \end{array} \right.
$$
and  $D^\alpha_0\circ \mathcal{C}_\beta^*(f)=\mathcal{C}_\beta^*\circ D^\alpha_0(f)$, where $D^\alpha_0 f(t)=\frac{1}{\Gamma(\alpha+1)}\vert
t\vert^\alpha W^\alpha_0f(t)$ for  $f\in {\mathcal S}$ and $t\in
\R$.

Note that we may write
$$
\mathcal{C}^*_\beta f(t)= \beta\int_{1}^\infty {(s-1)^{\beta-1}\over s^\beta}f(ts)ds,\,\, t\not=0,
$$
for $f\in \S$. The proof of next  result runs parallel to the proof of Theorem \ref{theorem3.6} and \ref{spectrr}.

\begin{theorem} \label{lemma2.2}
Let  $\alpha \ge 0$,  $\beta>0$, $1\le p<\infty$ and the  generalized dual Ces\'{a}ro operator $\mathcal{C}^*_\beta$ on $\mathcal{T}_p^{(\alpha)}(\vert t\vert^\alpha)$. Then

\begin{itemize}
\item[(i)] The operator $\mathcal{C}_\beta^*$ is bounded on $\mathcal{T}_p^{(\alpha)}(\vert t\vert^\alpha)$ and
$$
||\mathcal{C}_\beta^*||=\frac{\Gamma(\beta+1)\Gamma(1/p)}{\Gamma(\beta+1/p)}.
$$

\item[(ii)] The dual operator of $\mathcal{C}_\beta$ on $\mathcal{T}_p^{(\alpha)}(\vert t\vert^\alpha)$ is $\mathcal{C}_\beta^*$  on $\mathcal{T}_{p'}^{(\alpha)}(\vert t\vert^\alpha)$, i.e.
$$
\langle \mathcal{C}_\beta f,g \rangle_\alpha=\langle f,\mathcal{C}_\beta^*g \rangle_\alpha, \qquad f \in \mathcal{T}_p^{(\alpha)}(\vert t\vert ^\alpha), \quad  g \in\mathcal{T}_{p'}^{(\alpha)}(\vert t\vert^\alpha),
$$
where $\langle \quad,\quad \rangle_\alpha$ is given in
Proposition \ref{juanjo2} (iii).

\item[(iii)] If $f\in \mathcal{T}_p^{(\alpha)}(\vert t\vert ^\alpha)$, then
\begin{eqnarray}\label{integral3}
\mathcal{C}_\beta^* f(t)&=& \displaystyle\beta\int_{-\infty}^0 (e^{-r}-1)^{\beta-1}e^{-r(1-1/p-\beta)}T_{r,p}f(t)dr, \quad t\in \R,\end{eqnarray}
where the $C_0$-group $(T_{r, p})_{r\in \RR}$ is defined in Theorem \ref{groups}.

\item[(iv)]
$$
\sigma(\mathcal{C}_\beta^\ast)=\Gamma(\beta+1)\overline{\left\{{\Gamma({1\over p}+it)\over \Gamma(\beta+{1\over p}+it)} \ : \ t\in\R \right\}}.
$$
\end{itemize}
\end{theorem}

\begin{remark} {\rm Note that for $t=0$, by the integral representation
(\ref{integral3})
$$
\mathcal{C}_\beta^* f(0) = f(0) \displaystyle\beta\int^{\infty}_0 (1-e^{-r})^{\beta-1}dr=\infty, \quad f\in {\mathcal S}.
$$
}

\end{remark}

\section{Fourier transform and Ces\`{a}ro generalized operator}

\setcounter{theorem}{0}
\setcounter{equation}{0}

We remind the reader that the Fourier transform of a function $f$ in $L^1(\RR)$ is defined by
$$
\hat f(t):=\int_{-\infty}^{\infty}e^{-ixt} f(x)dx, \qquad t\in \R.
$$
It is well-known that $\hat f$ is continuous on $\RR$ and $\hat f(t)\to 0$ when $\vert t\vert\to \infty$ (the Riemann-Lebesgue lemma). In the case that $f\in L^p(\RR)$ for some $1<p\le 2$, the Fourier transform of $f$ is defined in terms of  a limit in the norm of $L^{p'}(\RR)$ of truncated integrals:
$$
\hat f:= \lim_{R\to \infty}\widehat{f\chi_{(-R, R)}}, \qquad \widehat{f\chi_{(-R, R)}}(t)=\int_{-R}^{R}e^{-ixt} f(x)dx, \qquad t\in \R,
$$
i.e., $\hat f \in L^{p'}(\RR)$ and
$
\lim_{R\to \infty}\Vert \hat f-\widehat{f\chi_{(-R, R)}}\Vert_{p'}=0
$
where ${1\over p}+{1\over p'}=1$ and $\chi_{(-R, R)}$ is the characteristic function of the interval $(-R, R)$, see for example  \cite[Vol 2, p.254]{[Zi]}. Then  the existence of $\hat f(t)$ is guaranteed only at almost every $t$ and $\hat f$ may be non continuous and the Riemann-Lebesgue lemma  could not hold (unlike the case when $f\in L^1(\R))$.

In case that $f\in L^p(\R)$ for some $2<p<\infty$, the  Fourier transform $\hat f$ cannot be defined as an ordinary function although $\hat f$ can be defined as a tempered distribution, see for example \cite[pp 19-30]{[SW]}.

In the next theorem, we consider the Fourier transform on the Sobolev space $\mathcal{T}_{p}^{(n)}(\vert t\vert^n)$.

\begin{theorem}
Take $1\le p\le 2$ and $n\in \NN$. Then  $\hat f \in \mathcal{T}_{p'}^{(n)}(\vert t\vert^n)$ for $f\in  \mathcal{T}_{p}^{(n)}(\vert t\vert^n)$ and ${1\over p}+{1\over p'}=1$.
\end{theorem}

\bgproof
Take $f\in  \mathcal{T}_{p}^{(n)}(\vert t\vert^n)$. Since $\mathcal{T}_{p}^{(n)}(\vert t\vert^n)\subset \mathcal{T}_{p}^{(j)}(\vert t\vert^j)$, we have that $x^jf^{(j)}\in L^p(\RR)$ for $0\le j \le n$. As
$$
(it)^n (\hat f)^{(n)}(t)=\sum_{j=0}^n (-1)^n{n\choose j}{n!\over j!}\widehat{x^jf^{(j)}}(t), \qquad \textrm{ $n\in\NN$, $t$ a.e. on $\R$,}
$$
(see for example \cite{[Zi]}), we conclude that  $(it)^n (\hat f)^{(n)} \in L^{p'}(\RR)$  and then $\hat f\in \mathcal{T}_{p'}^{(n)}(\vert t\vert^n)$.
\edproof

In what follows, we show that
$$
\widehat{{\mathcal C}_\beta(f)}={\mathcal C}_\beta^*(\widehat{f}), \qquad \hbox{ and }\qquad \widehat{{\mathcal C}_\beta^*(f)}={\mathcal C}_\beta(\widehat{f}), \qquad f\in L^p(\RR),
$$
for $1<p\le 2$ (Theorem \ref{conmutan}). This theorem extends the case $\beta=1$ formulated in \cite{Be} and proved in \cite{Mo02}.
Our approach looks like to be new and is based in the integral representations of ${\mathcal C}_\beta(f)$ and ${\mathcal
C}_\beta^*(f)$ given in Section 3.

\begin{lemma}\label{lemaconmu}
Let  $1\le p\le 2$ and  the family of operators  $(T_{t,p})_{t\in \RR}$ defined by
$T_{t,p}(f):=e^{-\frac{t}{p}}f(e^{-t}\cdot),$ for $f \in L^p(\RR)$. Then
$$
\widehat{T_{t,p} (f)}=T_{-t,p'}( \hat f), \qquad f  \in L^p(\RR), \qquad {1\over p}+{1\over p'}=1.
$$
\end{lemma}

\bgproof
Consider $1<p\le 2$,  $R>0$ and $f\in L^p(\R)$. We write by $f_R:=f\chi_{(-R, R)}$. Taking $\alpha=0$ in Theorem \ref{groups}, we conclude that $T_{t,p} (f)\in L^p(\R)$ and then $\widehat{T_{t,p} (f)}=\lim_{R\to \infty}\widehat{(T_{t,p} (f))_R}$ in $L^{p'}(\R)$. Note that
\begin{eqnarray*}
\widehat{(T_{t,p} (f))_R}(r)&=&e^{t\over p}\int_{-R}^{R}e^{-irx}f(e^{-t}x)dx= e^{t(1-{1\over p})}\int_{-e^{-t} R}^{e^{-t}R}e^{-i re^t y}f(y)dy= T_{-t, p'}(\widehat{f_{e^{-t}R}}).
\end{eqnarray*}
Since $\lim_{R\to \infty}\widehat{f_{e^{-t}R}}= \widehat{f}$ in $L^{p'}(\RR)$,  and $T_{t, p'}$ is continuous on $L^{p'}(\R)$, we conclude the result. The case $p=1$ may be considered directly taking $R=\infty$ in the integral expressions.
\edproof

\begin{remark}
{\rm Since $\mathcal{T}_p^{(\alpha)}(\vert t\vert^\alpha)\hookrightarrow  L^p(\mathbb{R})$ (Proposition \ref{juanjo2} (i)), the equality $\widehat{T_{t,p} (f)}=T_{t,p'}( \hat f)$ holds for $f\in\mathcal{T}_p^{(\alpha)}(\vert t\vert^\alpha)$ for $\alpha \ge 0$ and $1\le p\le 2$.
}
\end{remark}

Finally, we are ready to prove the main result in this section.

\begin{theorem}\label{conmutan}
Let $\beta >0$.
\begin{itemize}
\item[(i)] If $f\in L^p(\RR)$ for some $1<p\le 2$, then $\widehat{{\mathcal C}_\beta(f)}={\mathcal C}_\beta^*(\widehat{f})$.

\item[(ii)] If $f\in L^p(\RR)$ for some $1\le p\le 2$, then $\widehat{{\mathcal C}_\beta^*(f)}={\mathcal C}_\beta(\widehat{f})$.
\end{itemize}
\end{theorem}

\bgproof
(i) Take $f\in L^p(\RR)$ for some $1<p\le 2$. By Theorem \ref{lemma2.1} (ii) and Lemma \ref{lemaconmu} we have that
\begin{eqnarray*}
\widehat{{\mathcal C}_\beta(f)}(x)&=&\displaystyle\beta\int_0^\infty (1-e^{-r})^{\beta-1}e^{-r(1-1/p)}\widehat{T_{r,p}f}(x)dr\\
&=&\displaystyle\beta\int_{-\infty}^0 (e^{-r}-1)^{\beta-1}e^{-r(1/p-\beta)}T_{r,p'}\widehat{f}(x)dr\\
&=&\displaystyle\beta\int_{-\infty}^0 (e^{-r}-1)^{\beta-1}e^{-r(1-{1\over p'}-\beta)}T_{r,p'}\widehat{f}(x)dr={\mathcal C}_\beta^*(\widehat{f})(x)
\end{eqnarray*}
for almost every $x$ on $\R$ and we use the formula (\ref{integral3}).

(ii) Now take  $f\in L^p(\RR)$ for some $1\le p\le 2$. By the integral representation (\ref{integral3}) of ${\mathcal C}_\beta^*$ and  Lemma \ref{lemaconmu} we have that
\begin{eqnarray*}
\widehat{{\mathcal C}_\beta^*(f)}(x)&=&\displaystyle\beta\int_{-\infty}^0 (e^{-r}-1)^{\beta-1}e^{-r(1-{1\over p}-\beta)}T_{-r,p'}\widehat{f}(x)dr\\
&=&\beta\int_0^\infty (1-e^{-r})^{\beta-1}e^{-{r\over p}}T_{r,p'}\widehat{f}(x)dr\\
&=&\beta\int_0^\infty (1-e^{-r})^{\beta-1}e^{-r(1-{1\over p'})}T_{r,p'}\widehat{f}(x)dr= {\mathcal C}_\beta(\widehat{f})(x)\\
\end{eqnarray*}
for almost every $x$ on $\R$ and we use the Theorem \ref{lemma2.1} (ii).
\edproof

\begin{remark}{\rm
 By the Proposition \ref{continu}, we get that
$\widehat{{\mathcal C}_\beta(f)}(t)={\mathcal
C}_\beta^*(\widehat{f})(t)$ and  $ \widehat{{\mathcal
C}_\beta^*(f)}(t)={\mathcal C}_\beta(\widehat{f})(t)$ for
$t\not=0$ and $f \in \mathcal{T}_p^{(\alpha)}(\vert
t\vert^\alpha)$, $1<p\le 2$ and $\alpha \ge 1$. }
\end{remark}

\subsection*{Acknowledgements}

R. Ponce  wishes to thank the members of the Instituto
Universitario de Matem\'aticas y Aplicaciones (I.U.M.A.) at
Universidad de Zaragoza for their kind hospitality.

\end{document}